\documentclass[12pt]{article}
\usepackage{geometry}
\usepackage{amsfonts}
\usepackage{dsfont}
\usepackage[latin1]{inputenc}
\usepackage{amssymb , amsmath, amsthm, amsopn, amstext, amscd}
\usepackage[T1]{fontenc}
\usepackage{latexsym, enumerate}
\usepackage{graphics}
\usepackage{graphicx}

\allowdisplaybreaks

\newtheorem{definit}{Definition}[section]

\newtheorem{thm}{Theorem}
\newtheorem{lm}{Lemma}
\newtheorem{rmq}{Remark}

 \DeclareMathOperator*{\argmin}{Argmin}

\title{Generalization of $\ell_{1}$ constraints for high dimensional regression
problems}
\author{Pierre ALQUIER$^{(1)}$ and Mohamed HEBIRI$^{(2)}$ \\ \\
\small{(1, 2) LPMA, CNRS-UMR 7599,}\\ 
\small{Universit\'{e} Paris 7 - Diderot, UFR de Math\'{e}matiques,}\\ 
\small{175 rue de Chevaleret F-75013 Paris, France.}\\ \\
\small{(1) CREST-LS,}\\
\small{3, avenue Pierre Larousse}\\
\small{92240 Malakoff, France.}\\ \\
\small{(2) D\'epartement de Math\'ematiques, Universit\'e Paris-Est.}
}

\date{}

\begin{document}

\maketitle

\begin{abstract}
We focus on the high dimensional linear regression
$Y\sim\mathcal{N}(X\beta^{*},\sigma^{2}I_{n})$,
where $\beta^{*}\in\mathds{R}^{p}$ is the parameter of interest.
In this setting, several estimators such as the LASSO \cite{Tibshirani-LASSO}
and the Dantzig Selector \cite{Dantzig} are known to satisfy interesting
properties whenever the vector $\beta^{*}$ is sparse.
Interestingly both of the LASSO and the Dantzig Selector can be seen as
orthogonal projections of $0$
into $\mathcal{DC}(s)=\{\beta\in\mathds{R}^{p},\|X'(Y-X\beta)\|_{\infty}\leq
s\}$ - using an
$\ell_{1}$ distance for the Dantzig Selector and $\ell_{2}$ for the LASSO. For a
well chosen $s>0$, this set is actually a confidence region for $\beta^{*}$.
In this paper, we investigate the properties of estimators defined as
projections on $\mathcal{DC}(s)$ using general distances.
We prove that the obtained estimators satisfy oracle properties close to the one
of the LASSO and Dantzig Selector.
On top of that, it turns out that these estimators can be tuned to exploit a
different sparsity or/and slightly different estimation objectives.
\end{abstract}

\section{Introduction}\label{sec:Introd}

In many modern applications, one has to deal with very large datasets.
Regression problems may
involve a large number of covariates, possibly larger than the sample size. In
this situation,
a major issue lies in  dimension reduction which can be performed through the
selection of a small
amount of relevant covariates. For this purpose, numerous regression methods
have been proposed
in the literature, ranging from the classical information criteria such as
$\mathop{\rm C_{p}}$,
$\mathop{\rm AIC}$ and $\mathop{\rm BIC}$ to the more recent
regularization-based techniques such
as the $\ell_1$ penalized least square estimator, known as the
LASSO~\cite{Tibshirani-LASSO}, and the
Dantzig selector~\cite{Dantzig}. These $\ell_{1}$-regularized regression methods
have recently
witnessed several developments due to the attractive feature of  computational
feasibility, even
for high dimensional data when the number of covariates $p$ is large.

Consider the linear regression model
\begin{equation}
\label{model}
Y = X\beta^* + \varepsilon,
\end{equation}
where $Y$ is a vector in
$\mathbb{R}^n$, $\beta^*\in \mathbb{R}^p$ is the parameter vector, $X$ is an $n
\times p$
real-valued matrix with possibly much fewer rows than columns, $n \ll p$, and
$\varepsilon$ is
a random noise vector in $\mathbb{R}^n$. Here, for the sake of simplicity, we
will assume that
$\varepsilon\sim\mathcal{N}(0,\sigma^{2}I_{n})$. Let $\mathds{P}$ denote the
probability
distribution of $Y$ in this setting. Moreover, we assume that the matrix $X$ is
normalized in such
a way that $X'X$ has only $1$ on its diagonal.
The analysis of regularized regression methods
for high
dimensional data usually involves a sparsity assumption on $\beta^*$ through the
{\it sparsity
index} $\Vert \beta^*\Vert_0 = \sum_{j=1,\ldots,p} \mathbb{I}(\beta_j^* \neq 0)$
where
$\mathbb{I}(\cdot)$ is the indicator function. For any $q\geq 1$, $d\geq 0$ and
$a\in\mathbb{R}^{d}$, denote by $\|a\|_{q}^{q} = \sum_{i=1}^{d}|a_{i}|^{q}$ and
$\|a\|_{\infty}=\max_{1\leq i\leq d}|a_{i}|$, the $\ell_q$ and the $\ell_\infty$
norms respectively.
When the design matrix $X$ is normalized, the LASSO and the Dantzig selector
minimize respectively
$\Vert X \beta \Vert_2^2$ and $\Vert \beta \Vert_1 $ under the constraint
$ \Vert X'(Y-X\beta)\Vert_\infty \leq s$ where $s$ is a positive tuning
parameter
(e.g.~\cite{DualLasso,CSEL} for the dual form of the LASSO). This geometric
constraint is central
in the approach developed in the present paper and we shall use it in a general
perspective.
Let us mention that several objectives may be considered by the statistician
when we
deal with the model given by Equation~\eqref{model}.
Usually, we consider three specific objectives in the high-dimensional setting
(i.e.,
$p\ge n$):\\

\noindent {\bf Goal 1} - {\it Prediction:} The reconstruction of the signal $X
\beta^*$ with
the best possible accuracy is first considered. The quality of the
reconstruction with an
estimator $\hat{\beta}$ is often measured with the squared error $\| X
\hat{\beta} - X
\beta^*\|_2^2$. In the standard form, results are stated as follows: under
assumptions on
the matrix $X$ and with high probability, the prediction error is bounded by $ C
\log{(p)} \|
\beta^*\|_0 $ where $C$ is a positive constant. Such results for the prediction
issue have
been obtained in \cite{Lasso3,Bunea_consist,Lasso2} for the LASSO and in
\cite{Lasso3}
for the Dantzig selector. We also refer to
\cite{KoltchDant,Koltchl1plus,VanGpLass,
VandeGeerSparseLasso,ArnakTsyb,ChriMo7GpLass} for related works with different
estimators
(non-quadratic loss, penalties slightly different from $\ell_1$ and/or random
design). The
results obtained in the works above-mentioned are optimal up to a logarithmic
factor as it has
been proved in~\cite{BTWAggSOI}. See
also~\cite{VandeGeerConditionLasso09,cherno2} for very
nice survey papers, or the introduction of \cite{mohamoux}.\\

\noindent {\bf Goal 2} - {\it Estimation:} Another wishful thinking is that the
estimator
$\hat{\beta}$ is close to $\beta^*$ in terms of the $\ell_q$ distance for $q\geq
1$. The
estimation bound is of the form $ C \left\| \beta^*\right\|_0 (\log{(p)} /
n)^{q/2} $ where
$C$ is a positive constant. Such results are stated for the LASSO
in~\cite{BTWAggSOI,Lasso2}
when $q=1$, for the Dantzig selector in~\cite{Dantzig} when $q=2$ and have been
generalized
in~\cite{Lasso3} with $1\leq q \leq 2$ for both the LASSO end the Dantzig
selector.\\

\noindent {\bf Goal 3} - {\it Selection:} Since we consider variable selection
methods, the
identification of the true support $\{j:\,\beta_j^*\neq 0\}$ of the vector
$\beta^*$ is to
be considered. One expects that the estimator $\hat{\beta}$ and the true vector
$\beta^*$
share the same support at least when $n$ grows to infinity. This is known as the
variable
selection consistency problem and it has been considered for the LASSO and the
Dantzig Selector in several
works~\cite{Bunea_consist,KarimNormSup,MeinshBulhmConsistLasso,MeinYuSelect,
WainSelection,BiYuConsistLasso}.
\\

In this paper, we focus on variants of {\bf Goal 1} and {\bf Goal 2}, using
estimators $\hat{\beta}$ that
also satisfy the constraint $\Vert X'(Y-X\hat{\beta})\Vert_\infty \leq s$. It is
organized as
follows. In Section \ref{geometric} we give some general geometrical
considerations on the
LASSO and the Dantzig Selector that motivates the introduction of the general
form of estimator:
$$ \argmin_{\beta\in \Vert X'(Y-X\beta)\Vert_\infty \leq s} \|\beta \| $$
for any semi-norm $\|\cdot\|$.
In Section \ref{generalized}, we focus on two particular cases of interest in
this family,
and give some sparsity inequalities in the spirit of the ones in~\cite{Lasso3}.
We show that
under the hypothesis that $F\beta^{*}$ is sparse for a known matrix $F$, we are
able to
estimate properly $\beta^{*}$. Some application to a generic inverse problem are
provided
with numerical experiments. Finally, Section \ref{preuves} is dedicated to
proofs.

\section{Some geometrical considerations}
\label{geometric}

\begin{definit}
Let us put, for any $s>0$,
$ \mathcal{DC}(s) = \left\{\beta\in\mathds{R}^{p}: \Vert
X'(Y-X\beta)\Vert_\infty \leq s\right\}$.
\end{definit}

\begin{lm}
\label{lem}
For any $s>0$, $\mathds{P}(\beta^{*}\in\mathcal{DC}(s)) > 1-
p\exp(-s^{2}/(2\sigma^{2})) .$
\end{lm}

This means that $ \mathcal{DC}(s)$ is a confidence region for $\beta^{*}$.
Moreover, note that
$ \mathcal{DC}(s)$ is convex and closed. Let $\|\cdot\|$ be any semi-norm
in $\mathds{R}^{p}$.
Let $\Pi_{\|\cdot\|}^{s}$ denote an orthogonal projection on $ \mathcal{DC}(s)$
with respect to
$\|\cdot\|$:
$$ \Pi_{\|\cdot\|}^{s}(b) \in \argmin_{\beta\in \mathcal{DC}(s)} \|\beta-b\| .$$
From properties of projections, we know that
$$
 \beta^{*}\in\mathcal{DC}(s) \Rightarrow \forall b\in \mathds{R}^{p},
\|\Pi_{\|\cdot\|}^{s}(b)-\beta^*\|
 \leq \|b-\beta^{*}\| .$$

There is a very simple interpretation to this inequality: if $b$ is any
estimator of $\beta^{*}$,
then, with probability at least $1-p\exp(-s^{2}/(2\sigma^{2}))$,
$\Pi_{\|\cdot\|}^{s}(b)$ is a better
estimator.
In order to perform shrinkage it seems natural to take $b=0$.

\begin{definit}
We define our general estimator by
$$ \hat{\beta}^{\|\cdot\|}_{s} = \Pi_{\|\cdot\|}^{s}(0) \in \argmin_{\beta\in
\mathcal{DC}(s)} \|\beta\| .$$
\end{definit}

We have the following examples:
\begin{enumerate}
\item for $\|\cdot\|=\|\cdot\|_{1}$, we obtain the definition of the Dantzig
Selector given in \cite{Dantzig}.
\item for $\|\beta\|=\|X\beta\|_{2}$, we obtain the program
$ \argmin_{\beta\in \mathcal{DC}(s)} \|X\beta\|_{2}$.
It was proved in~\cite{DualLasso} for example that a particular solution of this
program is
Tibshirani's LASSO estimator~\cite{Tibshirani-LASSO} known as
$$ \hat{\beta}^{L}_{s} = \argmin_{\beta\in \mathds{R}^{p}} \left[
\|Y-X\beta\|_{2}^{2} + 2s\|\beta\|_{1}
\right].$$
\item for $\|\beta\|=\|X'X\beta\|_{q}$ with $q>0$, it is proved in \cite{CSEL}
that the
solution coincides with the "Correlation Selector" and it does not depend on
$q$.
\end{enumerate}

In the next Section, we exhibit other cases of interest and provide some
theoretical results
on the performances of the estimators.

\section{Generalized LASSO and Dantzig Selector}
\label{generalized}

\subsection{Definitions}

Let $F$ be an application $\mathds{R}_{+}\rightarrow\mathds{R}_{+}$
with the restriction that
$F(x)=0$ may be equal to $0$ only for $x=0$. Note that $X'X$ may be written, for
some orthogonal matrix $Q$,
$$ Q' \left(
\begin{array}{c c c}
 \lambda_1 & \dots & 0 \\ \vdots & \ddots & \vdots \\ 0 & \dots & \lambda_p
\end{array}
\right)Q,  \text{then we put } F(X'X) = Q' \left(
\begin{array}{c c c}
 F(\lambda_1) & \dots & 0 \\ \vdots & \ddots & \vdots \\ 0 & \dots &
F(\lambda_p)
\end{array}
\right)Q .$$
The idea is that, for a well chosen norm $\|\cdot\|$, we will build estimators that
will be useful
to estimate $\beta^{*}$ when $F(X'X)\beta^{*}$ is sparse, in the sense that they
will be close to
$\beta^{*}$ with respect to the semi-norm induced by $G(X'X)$ for $G(x)=xF(x)$.

\begin{definit}
We define the "Generalized Dantzig Selector", $ \hat{\beta}^{GDS}_{s} $, as $
\hat{\beta}^{\|\cdot\|}_{s} $ for
$\|b\|=\|F(X'X)b\|_{1}$, and the "Generalized LASSO", $ \hat{\beta}^{GL}_{s} $,
for
$\|b\|=(b'G(X'X)b)^{1/2} $.
\end{definit}

\begin{rmq}
\label{rmkcomplement}
In the case where the program
$ \min_{\beta\in \mathcal{DC}(s)} \beta'G(X'X)\beta $ has multiple solutions
we define $\hat{\beta}^{GL}_{s} $ as one of the solutions that minimizes
$\|F(X'X)\beta\|_{1}$
among all the solutions $\beta$. The case where the program
$ \min_{\beta\in \mathcal{DC}(s)} \|F(X'X)\beta\|_{1} $
has multiple solution does not cause any trouble: we can take $
\hat{\beta}^{GDS}_{s} $
as any of these solution without any effect on its statistical properties.
\end{rmq}

\subsection{Sparsity Inequalities}

We now present the assumptions we need to state the Sparsity Inequalities.

\noindent {\bf Assumption $A(c)$} for $c>0$:
for any $\alpha\in\mathds{R}^{p}$ such that
$$ \sum_{j:(F(X'X)\beta^{*})_{j}=0}
\left|\alpha_{j}\right| \leq 3 \sum_{j:(F(X'X)\beta^{*})_{j} \neq 0}
\left|\alpha_{j}\right|, $$ we have, for $H(x)=x/F(x)$ (with the convention
$0/0=0$),
\begin{equation*}
\sum_{j:(F(X'X)\beta^{*})_{j} \neq 0} \alpha_{j}^{2} \leq c \alpha' H(X'X)
\alpha.
\end{equation*}

This assumption can be seen as a modification of assumptions
in \cite{Lasso3}: if we put $F(x)=1$, $F(X'X)=I_p$ and $H(X'X)=X'X$ and we
obtain exactly the same
assumption that in \cite{Lasso3}.
For the sake of shorteness, we put $F=F(X'X)$, $G=G(X'X)$ and $H=H(X'X)$.

\begin{thm}
\label{legrostheoreme}
Let us take $\varepsilon\in]0,1[$ and $s=2\sigma(2\log(p/\varepsilon))^{1/2}$.
Assume that Assumption $A(c)$ is satisfied for some $c>0$.
With probability at least $1-\varepsilon$ we have simultaneously:
$$
\left\{
\begin{array}{r l}
(\hat{\beta}^{GDS}_{s}-\beta^{*})'G(\hat{\beta}^{GDS}_{s}-\beta^{*})
&
\leq
             72 \sigma^{2} c \|F\beta^{*}\|_{0}
\log(p/\varepsilon) ,
\\
 \|F(\hat{\beta}^{GDS}_{s}-\beta^{*})\|_{1}
&
\leq 18\sqrt{2}\sigma
\|F\beta^{*}\|_{0}
              \sqrt{c\log(p/\varepsilon)},
\\
(\hat{\beta}^{GL}_{s}-\beta^{*})'
G(\hat{\beta}^{GL}_{s}-\beta^{*})
&
\leq
128 \sigma^{2} c \|F\beta^{*}\|_{0} \log(p/\varepsilon),
\\
\|F(\hat{\beta}^{GL}_{s}-\beta^{*})\|_{1}
&
\leq 32\sqrt{2} \sigma
\|F\beta^{*}\|_{0}
              \sqrt{c\log(p/\varepsilon)}.
\end{array}
\right.
$$
\end{thm}

In the case $F(x)=1$, we obtain the same result as in \cite{Lasso3}.
However, it is worth noting that the use of $\hat{\beta}^{GL}_{s}$ is
particularly useful when $F\beta^{*}$ is sparse for a non-constant $F(x)$,
and $\beta^{*}$ is not.
In this case the errors of the LASSO and the Dantzig Selector are not controlled
anymore.
This generalization is also of some interests especially when
Assumption $A(c)$ is satisfied for $H$, but not satisfied if we
replace $H$ by $X'X$.
We now give an exemple.

\subsection{Application to a generic inverse problem}

In statistical inverse problems, one usually has to deal with the following
regression problem:
$ Y \sim \mathcal{N}(X\beta^{*},\sigma^{2} I_{n}) $ with a known $\sigma^{2}$,
$X$ a symmetric operator (for example a convolution operator)
and a regularity assumption on $\beta^{*}$. This assumption is often that
$\beta^{*}$ belongs to the
range of $X$ or of a power of $X$: $\beta^{*}= X^{\alpha}g$. See
for example \cite{cavalier}
and the references therein.

We will now assume that $g$ is sparse. In this case, note that $X'X=X^{2}$. As
$X^{-\alpha}\beta^*$ is sparse,
we put $F(x)=x^{-\alpha/2}$. So $F=X^{-\alpha}$ and $G=X^{2-\alpha}$.
In this case,
Theorem \ref{legrostheoreme} gives for example
$$ (\hat{\beta}^{GL}_{s}-\beta^{*})'G(\hat{\beta}^{GL}_{s}-\beta^{*}) \leq
             128 \sigma^{2} c \|g \|_{0} \log(p/\varepsilon),$$
under an assumption on $H=X^{2+\alpha}$ (it is worth mentionning
that in the case where
$\alpha=-2$, $H=I_n$ and so Assumption $A(c)$ is always satisfied with
$c=1$, even if the case $\alpha>0$ is more meaningful).

\begin{figure}[h]
\begin{center}
\includegraphics[height=4cm,width=14cm]{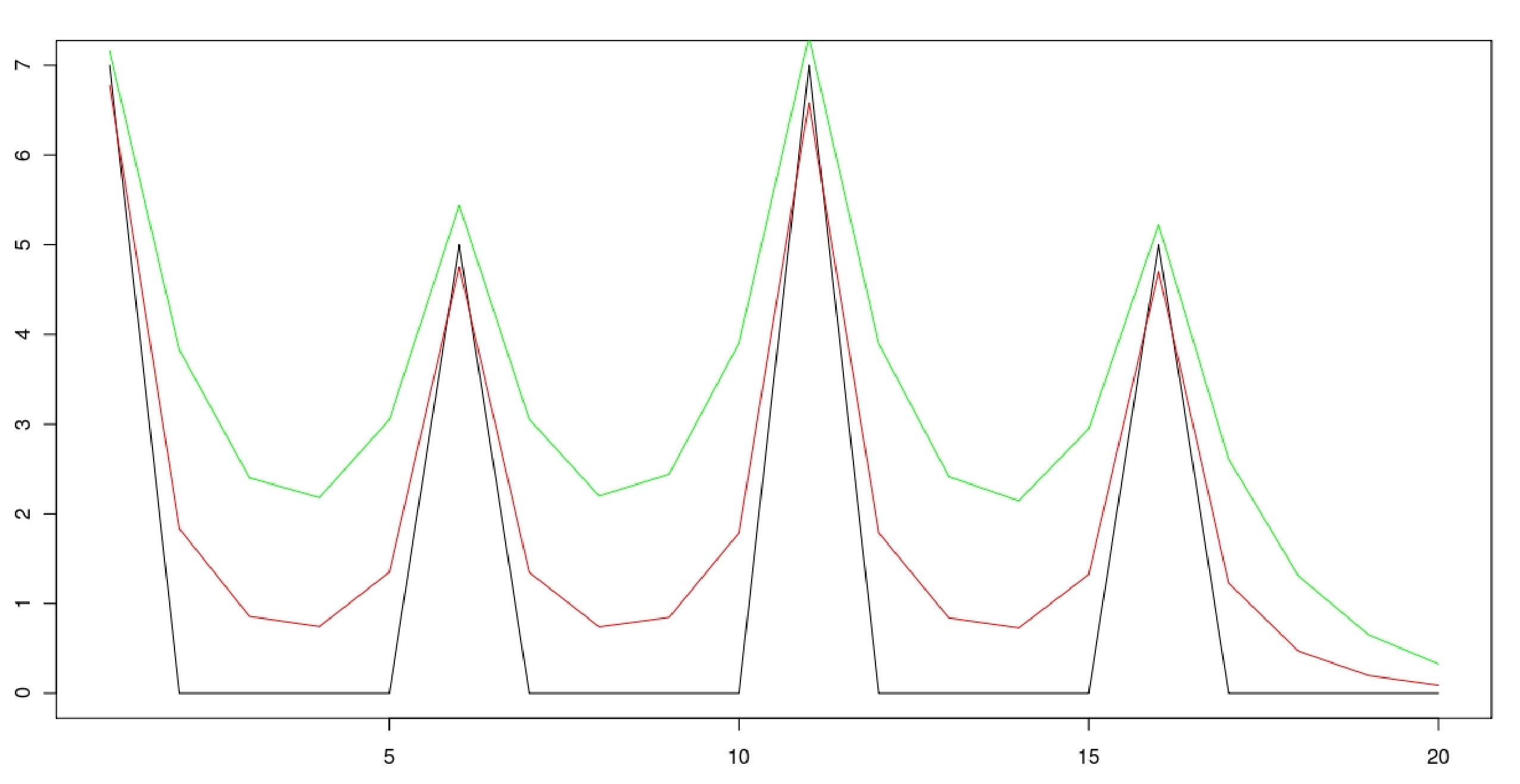}
\caption{The parameter $\beta^*=X^{\alpha}g$ for different values of $\alpha$. In black, $\alpha=0$,
so $\beta^* = g$ is sparse. In red, $\alpha=1$, $\beta^*$ is a bit smoothed, but still can be approximated
by a sparse signal. In green, $\alpha=2$, $\beta^*$ is smoother, and approximation by a sparse signal
do not hold any longer.}
\label{valalpha}
\end{center}
\end{figure}

We now provide a very short empirical comparison of the LASSO and Generalized
LASSO approach in a toy example of such a model. Note that for $\alpha\geq 0$,
$\beta^* = X^{\alpha}g$ being a smoothed version of $g$, is ``almost sparse'',
so a comparison with the LASSO makes sense. We propose the following setting:
let $M(\rho)=(\rho^{|i-j|})_{1\leq i,j\leq n}$, and
$X=M^{1/2}$.  We take
$g=(7,0,0,0,0,5,0,0,0,0,7,0,0,0,0,5,0,0,0,0)$, $n=20$, $\rho=0.5$, Figure \ref{valalpha} gives the different
values of $\beta^*=X^{\alpha}g$ for various values of $\alpha$.

We compute the LASSO and Generalized LASSO in each case, and report the
performance of the oracle with respect to the regularization parameter $s$:
$$ {\rm Perf.GL} = \inf_{s>0} \|X\hat{\beta}^{GL}_{s} - X\beta^*\| \text{ and }
 {\rm Perf.L} = \inf_{s>0} \|X\hat{\beta}^{L}_{s} - X\beta^*\| .$$
Of course, in practice, the optimal $s$ in unknown and may be estimated
by cross-validation for example.
We test both estimators with several values for the parameters
$\alpha$ and $\sigma^{2}$. For each value of these parameters, we run
$20$ experiments and report the mean performances for both estimators.
The results are given in Table \ref{lesresultats}. We can see that the results
seem coherent with Theorem \ref{legrostheoreme}: there seems to be an
advantage in practice to consider the Generalized LASSO in the cases
where $alpha \neq 0$.

\begin{table}[t!]
\caption{The mean results for $20$ experiments for each value of $(\alpha,\sigma^{2})$.}
\label{lesresultats}
\begin{center}
\begin{tiny}
\begin{tabular}{|p{1.0cm}|p{1.0cm}||p{2.0cm}|p{2.0cm}|}
\hline
$\alpha$  & $\sigma^{2}$ & mean of ${\rm Perf.L}$ & mean of ${\rm Perf.GL}$ \\
\hline \hline
                   & 0.01  & 0.167 & 0.118 \\
   -2              & 0.30  & 4.792 & 3.076 \\
                   & 1.00  & 16.636& 10.328\\ 
\hline
                   & 0.01  & 0.194 & 0.097 \\
   -1              & 0.30  & 5.624 & 2.911 \\
                   & 1.00  & 14.386& 8.56  \\
\hline
                   & 0.01  & 0.098 & 0.098 \\
  0                & 0.30  & 2.835 & 2.835 \\
                   & 1.00  & 9.012 & 9.012 \\ 
\hline
                   & 0.01  & 0.196 & 0.094 \\
   +1              & 0.30  & 5.144 & 2.517 \\
                   & 1.00  & 13.232& 8.597 \\ 
\hline
                   & 0.01  & 0.199 & 0.101 \\
   +2              & 0.30  & 5.589 & 3.018 \\
                   & 1.00  & 17.957& 10.228\\ 
\hline
                   & 0.01  & 0.183 & 0.102 \\
   +3              & 0.30  & 5.538 & 3.175 \\
                   & 1.00  & 19.133& 10.371\\ 
\hline
\end{tabular}
\end{tiny}
\end{center}
\vspace*{-6pt}
\end{table}

\section{Proofs}
\label{preuves}

\subsection{Proof of Lemma \ref{lem}}

We have $Y\sim\mathcal{N}(X\beta^{*},\sigma^{2}I_{n})$ and so
$Y-X\beta^{*}\sim\mathcal{N}(0,\sigma^{2}I_{n})$
and finally $X'(Y-X\beta^{*})\sim\mathcal{N}(0,\sigma^{2}X'X)$. Let us put
$V=X'(Y-X\beta^{*})$ and let $V_{j}$
denote the $j$-th coordinate of $V$. Note that $X'X$ is normalized such that
for any $j$, $V_{j}\sim\mathcal{N}(0,\sigma^{2})$, so:
$\mathds{P}\left(|V_{j}|>s\right) \leq \exp(-s^{2}/(2\sigma^{2}))$.
Then
$ \mathds{P}\left(\|V\|_{\infty}>s\right) \leq p \exp(-s^{2}/(2\sigma^{2}))$.
\qed

\subsection{Proof of Theorem \ref{legrostheoreme}}

We use arguments from \cite{Lasso3}.
From now, we assume that the event $\{\beta^{*}\in\mathcal{DC}(s/2)\} =
\{\|X'(Y-X\beta^{*}\|_{\infty}<s/2\}$ is
satisfied. According to Lemma \ref{lem}, the probability of this event is at
least
$1-p\exp(-s^{2}/(8\sigma^{2}))=1-\varepsilon$ as
$s=2(2\log(p/\varepsilon))^{1/2}$.

\noindent {\bf Proof of the results on the Generalized Dantzig Selector.}

We have
\begin{multline*}
(\hat{\beta}^{GDS}_{s}-\beta^{*})'G(\hat{\beta}^{GDS}_{s}-\beta^{*})
= (\hat{\beta}^{GDS}_{s}-\beta^{*})'X'XF(\hat{\beta}^{GDS}_{s}-\beta^{*})
\\
\leq
\|X'X(\hat{\beta}^{GDS}_{s}-\beta^{*})\|_{\infty}\|F(\hat{\beta}^{GDS}_{s}
-\beta^{*})\|_{1}
\\
\leq
\left(\|X'(Y-X\beta^{*})\|_{\infty} +
\|X'(Y-X\hat{\beta}^{GDS}_{s})\|_{\infty}\right)
\|F(\hat{\beta}^{GDS}_{s}-\beta^{*})\|_{1}
\\
\leq
(s/2+s) \|F(\hat{\beta}^{GDS}_{s}-\beta^{*})\|_{1}
\end{multline*}
since $\hat{\beta}^{GDS}_{s}\in\mathcal{DC}(s)$, and
$\{\beta^{*}\in\mathcal{DC}(s/2)\}$ is satisfied.
By definition of $\hat{\beta}^{GDS}_{s}$,
\begin{multline*}
0 \leq \|F\beta^{*}\|_{1} -\|F\hat{\beta}^{GDS}_{s}\|_{1}
\\
       = \sum_{(F\beta^{*})_{j}\neq 0} |(F\beta^{*})_{j}| -
\sum_{(F\beta^{*})_{j}\neq 0} |(F\hat{\beta}^{GDS}_{s})_{j}|
                      - \sum_{(F\beta^{*})_{j} = 0}
|(F\hat{\beta}^{GDS}_{s})_{j}|
\\
\leq \sum_{(F\beta^{*})_{j}\neq 0} |(F\beta^{*})_{j}
-(F\hat{\beta}^{GDS}_{s})_{j}|
           - \sum_{(F\beta^{*})_{j}= 0} |(F\beta^{*})_{j}
-(F\hat{\beta}^{GDS}_{s})_{j}|.
\end{multline*}
This means that
$$
\|F(\hat{\beta}^{GDS}_{s}-\beta^{*})\|_{1} \leq 2 \sum_{(F\beta^{*})_{j}\neq 0}
     |(F\beta^{*})_{j} -(F\hat{\beta}^{GDS}_{s})_{j}|.
$$
We can summarize all that we have now:
\begin{multline}
\label{etapepremiere}
(\hat{\beta}^{GDS}_{s}-\beta^{*})'G(\hat{\beta}^{GDS}_{s}-\beta^{*})
\leq \frac{3s}{2} \|F(\hat{\beta}^{GDS}_{s}-\beta^{*})\|_{1}
\\
\leq 3s \sum_{(F\beta^{*})_{j}\neq 0}
     |(F\beta^{*})_{j} -(F\hat{\beta}^{GDS}_{s})_{j}|.
\end{multline}
Let us remark that Inequality \eqref{etapepremiere} implies that the vector
$\alpha=F(\hat{\beta}^{GDS}_{s}-\beta^{*})$ may be used in Assumption $A(c)$.
This leads to
\begin{multline}
\label{etapedeuxieme}
(\hat{\beta}^{GDS}_{s}-\beta^{*})'G(\hat{\beta}^{GDS}_{s}-\beta^{*})
\leq
3s \sum_{(F\beta^{*})_{j}\neq 0}
     |(F\beta^{*})_{j} -(F\hat{\beta}^{GDS}_{s})_{j}|
\\
\leq 3s \sqrt{\|F\beta^{*}\|_{0}\sum_{(F\beta^{*})_{j}\neq 0}
     [(F\beta^{*})_{j} -(F\hat{\beta}^{GDS}_{s})_{j}]^{2} }
\\
\leq 3s \sqrt{\|F\beta^{*}\|_{0} c
(F\hat{\beta}^{GDS}_{s}-F\beta^{*})'H(F\hat{\beta}^{GDS}_{s}-F\beta^{*})}
\\
= 3s \sqrt{\|F\beta^{*}\|_{0} c (\hat{\beta}^{GDS}_{s}-\beta^{*})'
G(\hat{\beta}^{GDS}_{s}-\beta^{*})}.
\end{multline}
As a consequence,
$$ (\hat{\beta}^{GDS}_{s}-\beta^{*})'A(\hat{\beta}^{GDS}_{s}-\beta^{*}) \leq 9
s^{2} \|P\beta^{*}\|_{0} c
                        = 72 \sigma^{2} c \|P\beta^{*}\|_{0}
\log(p/\varepsilon) .$$
Plugging this result into Inequality \eqref{etapedeuxieme} and using Inequality
\eqref{etapepremiere}
again, we obtain:
$$ \|P(\hat{\beta}^{GDS}_{s}-\beta^{*})\|_{1} \leq 18\sqrt{2}\sigma
\|P\beta^{*}\|_{0}
              \sqrt{c\log(p/\varepsilon)} .$$

\noindent {\bf Proof of the results on the Generalized LASSO.} \\
\noindent {\bf Step 1.}
As a fist step, we establish an important property of the Generalized LASSO estimator. We
prove that
\begin{multline}
\label{formedualeLASSO}
\forall \beta\in\mathds{R}^{p},\quad
|Y-XF\hat{\beta}^{GL}_{s}\|_{2}^{2} + 2s\|F\hat{\beta}^{GL}_{s}\|_{1} +
(\hat{\beta}^{GL}_{s})'F(H-X'X)F
\hat{\beta}^{GL}_{s}
\\
\leq
\|Y-XF\beta\|_{2}^{2} + 2s\|F\beta\|_{1} + \beta'F(H-X'X)F\beta.
\end{multline}
To prove Inequality \eqref{formedualeLASSO}, we write the Lagrangian of the
program that defines $\hat{\beta}^{GL}_{s}$:
\begin{equation*}
\mathcal{L}(\beta,\lambda,\mu) = \beta'G\beta
     + \lambda'  \left[ X'(X\beta-Y) - s E \right]
     + \mu' \left[ X'(Y-X\beta) - s E \right],
\end{equation*}
where $E=(1,\ldots,1)'$, $\lambda$ and $\mu$ are vectors in $\mathds{R}^{p}$.
Any solution $\underline{\beta}=\underline{\beta}(\lambda,\mu)$ must satisfy,
for some
$\lambda_{j}\geq 0$, $\mu_{j} \geq 0$ and $\lambda_{j}\mu_{j}=0$,
$$
0 = \frac{\partial \mathcal{L}}{\partial \beta}(\underline{\beta},\lambda,\mu)
= 2G\beta + X'X(\lambda-\mu),
$$
and then $
 G\underline{\beta} = (X'X)(\mu-\lambda)/2$ .
Note that $\lambda_{j}\geq 0$, $\mu_{j} \geq 0$ and $\lambda_{j}\mu_{j}=0$
imply that there is a $\gamma_{j}\in\mathds{R}$ such that
$\gamma_{j}=(\mu_{j}-\lambda_{j})/2$,
$|\gamma_{j}| = (\lambda_{j}+\mu_{j})/2$. Hence $\lambda_{j}=2(\gamma_{j})_{-}$
and $\mu_{j}=2(\gamma_{j})_{+}$, where for any $a$, $(a)_{+}=max(a;0)$ and 
$(a)_{-}=max(-a;0)$.
Let also $\gamma$ denote the vector which $j$-th component
is exactly $\gamma_{j}$, we obtain:
\begin{equation}
\label{SOLU}
 G\underline{\beta} = (X'X)\gamma.
\end{equation}
Then we have easily
$
\underline{\beta}'F\underline{\beta}
 = \underline{\beta}'(X'X)\gamma
 = \gamma'H\gamma
$.
Using these relations, the Lagrangian may be written:
\begin{eqnarray*}
\mathcal{L}(\underline{\beta},\lambda,\mu) &
= &  \gamma'H\gamma
 + 2 \gamma' X' Y - 2\gamma' (X'X) \underline{\beta}
- 2 s \sum_{j=1}^{p} |\gamma_{j}|
\\
& = & 2 \gamma' X' Y - \gamma 'H \gamma
 - 2 s \left\|\gamma\right\|_{1}
\end{eqnarray*}
Note that $\lambda$ and $\beta$, and so $\gamma$, should maximize this value.
Hence, $\gamma$ is to minimize
$$ - 2 \gamma' X' Y + \gamma 'H \gamma
 + 2 s \|\gamma\|_{1} + Y'Y
$$
Now, note that
$$
Y'Y - 2 \gamma' X' Y = \|Y-X\gamma\|_{2}^{2} -\gamma'(X'X)\gamma
$$
and then $\gamma$ also minimizes
$$
\|Y-X\gamma\|_{2}^{2} + 2 s \left\|\gamma\right\|_{1} + \gamma'
\left[H-(X'X)\right]\gamma.
$$
We end the proof of~\eqref{formedualeLASSO} by noting that
for every $b$ such that $Fb=\gamma$, then $b$ is to minimize
\begin{equation}
\label{progb}
\|Y-XFb\|_{2}^{2} + 2 s \left\|Fb\right\|_{1} + (Fb)'
\left[H-(X'X)\right](Fb).
\end{equation}
and that $\hat{\beta}^{GL}_{s}$ is such a $b$.

\noindent {\bf Step 2.} The next step is to apply Equation \eqref{formedualeLASSO} with
$\beta=\beta^{*}$ to obtain
\begin{multline*}
\|Y-XF\hat{\beta}^{GL}_{s}\|_{2}^{2} + 2s\|F\hat{\beta}^{GL}_{s}\|_{1} +
(\hat{\beta}^{GL}_{s})'F(H-X'X)F
\hat{\beta}^{GL}_{s}
\\
\leq
\|Y-XF\beta^{*}\|_{2}^{2} + 2s\|F\beta^{*}\|_{1} +
(F\beta^{*})'(H-X'X)F\beta^{*}.
\end{multline*}
For the sake of simplicity, we can define $\hat{\gamma}=F\hat{\beta}^{GL}_{s}$
(following the notations of Step 1)
and $\gamma^{*}=F\beta^{*}$ and
we obtain
\begin{multline*}
\|Y-X\hat{\gamma}\|_{2}^{2} + 2s\|\hat{\gamma}\|_{1} +
\hat{\gamma}'(H-X'X)\gamma
\\
\leq
\|Y-X\gamma^{*}\|_{2}^{2} + 2s\|\gamma^{*}\|_{1} +
(\gamma^{*})'(H-X'X)\gamma^{*}.
\end{multline*}
Computations lead to
\begin{multline*}
\|X(\hat{\gamma}-\gamma^{*})\|_{2}^{2} + 2s\|\hat{\gamma}\|_{1} +
\hat{\gamma}'(H-X'X)\hat{\gamma}
 -2(Y-X\gamma^{*})'X\hat{\gamma}
\\
+2(\gamma^{*})'(H-X'X)(\gamma^{*}-\gamma)
\leq
2s\|\gamma^{*}\|_{1} + (\gamma^{*})'(H-X'X)\hat{\gamma}
 -2(Y-X\gamma^{*})'X\gamma^{*},
\end{multline*}
and then
\begin{multline*}
\|X(\hat{\gamma}-\gamma^{*})\|_{2}^{2}
\\
\leq 2s (\|\gamma^{*}\|_{1}-\|\hat{\gamma}\|_{1})
           + 2(Y-X\gamma^{*})'X(\hat{\gamma}-\gamma^{*})
              -
(\gamma^{*}-\hat{\gamma})'(H-X'X)(\gamma^{*}-\hat{\gamma}).
\end{multline*}
As a consequence
\begin{multline*}
(\gamma^{*}-\hat{\gamma})'H(\gamma^{*}-\hat{\gamma})
\leq
2s (\|\gamma^{*}\|_{1}-\|\hat{\gamma}\|_{1})
           + 2(Y-X\gamma^{*})'X(\hat{\gamma}-\gamma^{*})
\\
\leq 2s\sum_{j=1}^{p} (|\gamma_{j}^{*}|-|\hat{\gamma}_{j}|) + 2
\|X'(Y-X\beta^{*})\|_{\infty}\sum_{j=1}^{p}
                           |\hat{\gamma}_{j}-\gamma^{*}_{j}|
\\
\leq 2s\sum_{j=1}^{p} (|\gamma_{j}^{*}|-|\hat{\gamma}_{j}|) + s \sum_{j=1}^{p}
                           |\hat{\gamma}_{j}-\gamma^{*}_{j}|.
\end{multline*}
So we obtain
\begin{multline}
\label{pivot}
(\gamma^{*}-\hat{\gamma})'H(\gamma^{*}-\hat{\gamma})
+ s\sum_{j=1}^{p} |\hat{\gamma}_{j}-\gamma^{*}_{j}|
\leq
2s\sum_{j=1}^{p} (|\hat{\gamma}_{j}|-|\gamma^{*}_{j}|) + 2 s \sum_{j=1}^{p}
                           |\hat{\gamma}_{j}-\gamma^{*}_{j}|
\\
= 2s \sum_{j:\gamma^{*}_{j}\neq 0} (|\hat{\gamma}_{j}|-|\gamma^{*}_{j}|) + 2 s
\sum_{j:\gamma^{*}_{j}\neq 0}
                           |\hat{\gamma}_{j}-\gamma^{*}_{j}|
= 4s \sum_{j:\gamma^{*}_{j}\neq 0}  |\hat{\gamma}_{j}-\gamma^{*}_{j}|.
\end{multline}
In particular, Equation \eqref{pivot} implies that
$$ \sum_{j:\gamma^{*}_{j}=0} |\hat{\gamma}_{j}-\gamma^{*}_{j}| \leq
           3 \sum_{j:\gamma^{*}_{j}\neq 0}  |\hat{\gamma}_{j}-\gamma^{*}_{j}|,$$
and so $\alpha=\hat{\gamma}_{j}-\gamma^{*}_{j}$ may be used in Assumption
$A(c)$. Then
Inequality \eqref{pivot} becomes
\begin{multline*}
(\gamma^{*}-\hat{\gamma})'H(\gamma^{*}-\hat{\gamma})
\leq
4s \sum_{j:\gamma^{*}_{j}\neq 0}  |\hat{\gamma}_{j}-\gamma^{*}_{j}|
\leq
4s \sqrt{
\|\gamma^{*}\|_{0} \sum_{j:\gamma^{*}_{j}\neq 0} 
(\hat{\gamma}_{j}-\gamma^{*}_{j})^{2} }
\\
\leq
4s \sqrt{
\|\gamma^{*}\|_{0} c
(\gamma^{*}-\hat{\gamma})'
H(\gamma^{*}-\hat{\gamma})}.
\end{multline*}
That leads to
\begin{equation}
\label{resultatpresquefinal}
(\hat{\beta}^{GL}_{s}-\beta^{*})'
G(\hat{\beta}^{GL}_{s}-\beta^{*})
=
(\gamma^{*}-\hat{\gamma})'H(\gamma^{*}-\hat{\gamma})
\leq
128 \sigma^{2} c \|F\beta^{*}\|_{0} \log(p/\varepsilon).
\end{equation}
We plug \eqref{resultatpresquefinal} into \eqref{pivot} again to obtain
$
\|\hat{\gamma}-\gamma^{*}\|_{1} \leq 32\sqrt{2} \sigma \|P\beta^{*}\|_{0}
              \sqrt{c\log(p/\varepsilon)}$.
\qed


\bibliographystyle{alpha}      
\renewcommand{\refname}{}
\bibliography{GenLasDan2}   

\begin{thebibliography}{MVdGB09}

\bibitem[Alq08]{CSEL}
P.~Alquier.
\newblock Lasso, iterative feature selection and the correlation selector:
  Oracle inequalities and numerical performances.
\newblock {\em Electron. J. Stat.}, pages 1129--1152, 2008.

\bibitem[BC11]{cherno2}
A.~Belloni and V.~Chernozhukov.
\newblock High dimensional sparse econometric models: An introduction.
\newblock In P.~Alquier, E.~Gautier, and G.~Stoltz, editors, {\em Inverse
  Problems and High-Dimensional Estimation}. Springer Lecture Notes in
  Statistics, 2011.

\bibitem[BRT09]{Lasso3}
P.~Bickel, Y.~Ritov, and A.~Tsybakov.
\newblock Simultaneous analysis of lasso and {D}antzig selector.
\newblock {\em Ann. Statist.}, 37(4):1705--1732, 2009.

\bibitem[BTW07a]{BTWAggSOI}
F.~Bunea, A.~Tsybakov, and M.~Wegkamp.
\newblock Aggregation for {G}aussian regression.
\newblock {\em Ann. Statist.}, 35(4):1674--1697, 2007.

\bibitem[BTW07b]{Lasso2}
F.~Bunea, A.~Tsybakov, and M.~Wegkamp.
\newblock Sparsity oracle inequalities for the lasso.
\newblock {\em Electron. J. Stat.}, 1:169--194, 2007.

\bibitem[Bun08]{Bunea_consist}
F.~Bunea.
\newblock {\em Consistent selection via the Lasso for high dimensional
  approximating regression models}, volume~3.
\newblock IMS Collections, 2008.

\bibitem[Cav11]{cavalier}
L.~Cavalier.
\newblock Inverse problems in statistics.
\newblock In P.~Alquier, E.~Gautier, and G.~Stoltz, editors, {\em Inverse
  Problems and High-Dimensional Estimation}. Springer Lecture Notes in
  Statistics, 2011.

\bibitem[CH08]{ChriMo7GpLass}
C.~Chesneau and M.~Hebiri.
\newblock Some theoretical results on the grouped variables lasso.
\newblock {\em Mathematical Methods of Statistics}, 17(4):317--326, 2008.

\bibitem[CT07]{Dantzig}
E.~Candes and T.~Tao.
\newblock The dantzig selector: statistical estimation when $p$ is much larger
  than $n$.
\newblock {\em Ann. Statist.}, 35, 2007.

\bibitem[DT07]{ArnakTsyb}
A.~Dalalyan and A.B. Tsybakov.
\newblock Aggregation by exponential weighting and sharp oracle inequalities.
\newblock {\em COLT 2007 Proceedings. Lecture Notes in Computer Science 4539
  Springer}, pages 97--111, 2007.

\bibitem[Heb09]{mohamoux}
M.~Hebiri.
\newblock {\em Quelques questions de s\'election de variables autour de
  l'estimateur {LASSO}}.
\newblock PhD thesis, 2009.

\bibitem[Kol09a]{KoltchDant}
V.~Koltchinskii.
\newblock The {D}antzig selector and sparsity oracle inequalities.
\newblock {\em Bernoulli}, 15(3):799--828, 2009.

\bibitem[Kol09b]{Koltchl1plus}
V.~Koltchinskii.
\newblock Sparse recovery in convex hulls via entropy penalization.
\newblock {\em Ann. Statist.}, 37(3):1332--1359, 2009.

\bibitem[Lou08]{KarimNormSup}
K.~Lounici.
\newblock Sup-norm convergence rate and sign concentration property of {L}asso
  and {D}antzig estimators.
\newblock {\em Electron. J. Stat.}, 2:90--102, 2008.

\bibitem[MB06]{MeinshBulhmConsistLasso}
N.~Meinshausen and P.~B{\"u}hlmann.
\newblock High-dimensional graphs and variable selection with the lasso.
\newblock {\em Ann. Statist.}, 34(3):1436--1462, 2006.

\bibitem[MVdGB09]{VanGpLass}
L.~Meier, S.~Van~de Geer, and P.~B{\"u}hlmann.
\newblock High-dimensional additive modeling.
\newblock {\em Ann. Statist.}, 37(6B):3779--3821, 2009.

\bibitem[MY09]{MeinYuSelect}
N.~Meinshausen and B.~Yu.
\newblock Lasso-type recovery of sparse representations for high-dimensional
  data.
\newblock {\em Ann. Statist.}, 37(1):246--270, 2009.

\bibitem[OPT00]{DualLasso}
M.~Osborne, B.~Presnell, and B.~Turlach.
\newblock On the {LASSO} and its dual.
\newblock {\em J. Comput. Graph. Statist.}, 9(2):319--337, 2000.

\bibitem[Tib96]{Tibshirani-LASSO}
R.~Tibshirani.
\newblock Regression shrinkage and selection via the lasso.
\newblock {\em J. Roy. Statist. Soc. Ser. B}, 58(1):267--288, 1996.

\bibitem[vdG08]{VandeGeerSparseLasso}
S.~van~de Geer.
\newblock High-dimensional generalized linear models and the lasso.
\newblock {\em Ann. Statist.}, 36(2):614--645, 2008.

\bibitem[vdGB09]{VandeGeerConditionLasso09}
S.~van~de Geer and P.~B{\"u}hlmann.
\newblock On the conditions used to prove oracle results for the lasso.
\newblock {\em Elect. Journ. Statist.}, 3:1360--1392, 2009.

\bibitem[Wai06]{WainSelection}
M.~Wainwright.
\newblock Sharp thresholds for noisy and high-dimensional recovery of sparsity
  using l1-constrained quadratic programming.
\newblock Technical report n. 709, Department of Statistics, UC Berkeley, 2006.

\bibitem[ZY06]{BiYuConsistLasso}
P.~Zhao and B.~Yu.
\newblock On model selection consistency of {L}asso.
\newblock {\em J. Mach. Learn. Res.}, 7:2541--2563, 2006.

\end{thebibliography}

\end{document}